\documentclass[12pt,oneside,english]{amsart}
\usepackage[T1]{fontenc}
\usepackage[latin9]{inputenc}
\usepackage[letterpaper]{geometry}
\usepackage{amsthm}
\usepackage{setspace}
\usepackage{amssymb}

\makeatletter

\newcommand{\noun}[1]{\textsc{#1}}

\numberwithin{equation}{section} 
\numberwithin{figure}{section} 

\makeatother

\usepackage{babel}

\begin{document}

\title{\textit{\noun{Effective $H^{\infty}$ interpolation constrained by
weighted Hardy and Bergman norms}}}

\author{Rachid Zarouf}
\begin{abstract}
Given a finite subset $\sigma$ of the unit disc $\mathbb{D}$ and
a holomorphic function $f$ in $\mathbb{D}$ belonging to a class
$X$, we are looking for a function $g$ in another class $Y$ which
satisfies $g_{\vert\sigma}=f_{\vert\sigma}$ and is of minimal norm
in $Y$. More precisely, we consider the interpolation constant $c\left(\sigma,\, X,\, Y\right)=\mbox{sup}{}_{f\in X,\,\parallel f\parallel_{X}\leq1}\mbox{inf}_{g_{\vert\sigma}=f_{\vert\sigma}}\left\Vert g\right\Vert _{Y}.$
When $Y=H^{\infty}$, our interpolation problem includes those of
Nevanlinna-Pick (1916) and Caratheodory-Schur (1908). If $X$ is a
Hilbert space belonging to the families of weighted Hardy and Bergman
spaces, we obtain a sharp upper bound for the constant $c\left(\sigma,\, X,\, H^{\infty}\right)$
in terms of $n=\mbox{card}\,\sigma$ and $r=\mbox{max}{}_{\lambda\in\sigma}\left|\lambda\right|<1$.
If $X$ is a general Hardy-Sobolev space or a general weighted Bergman
space (not necessarily of Hilbert type), we also establish upper and
lower bounds for $c\left(\sigma,\, X,\, H^{\infty}\right)$ but with
some gaps between these bounds. This problem of constrained interpolation
is partially motivated by applications in matrix analysis and in operator
theory.
\end{abstract}
\maketitle

\section*{1. Introduction}

\subsection*{a. Statement and historical context of the problem}

Let $\mathbb{D}=\{z\in\mathbb{C}:\,\vert z\vert<1\}$ be the unit
disc of the complex plane and let ${\rm Hol}\left(\mathbb{D}\right)$
be the space of holomorphic functions on $\mathbb{D}.$ We consider
here the following problem: given two Banach spaces $X$ and $Y$
of holomorphic functions on the unit disc $\mathbb{D},$ $X,\, Y\subset{\rm Hol}\left(\mathbb{D}\right),$
and a finite subset $\sigma$ of $\mathbb{D}$, what is the best possible
interpolation by functions of the space $Y$ for the traces $f_{\vert\sigma}$
of functions of the space $X$, in the worst case? The case $X\subset Y$
is of no interest, and so one can suppose that either $Y\subset X$
or $X$ and $Y$ are incomparable. Here and later on, $H^{\infty}$
stands for the space (algebra) of bounded holomorphic functions in
the unit disc $\mathbb{D}$ endowed with the norm $\left\Vert f\right\Vert _{\infty}=\sup_{z\in\mathbb{D}}\left|f(z)\right|.$ 

More precisely, our problem is to compute or estimate the following
interpolation constant\[
c\left(\sigma,\, X,\, Y\right)={\displaystyle \sup_{f\in X,\,\parallel f\parallel_{X}\leq1}}\mbox{inf}\left\{ \left\Vert g\right\Vert _{Y}:\, g_{\vert\sigma}=f_{\vert\sigma}\right\} \,.\]
For $r\in[0,\,1)$ and $n\geq1,$ we also define

\vspace{0.03cm}
\[
C_{n,\, r}(X,Y)=\mbox{sup}\left\{ c(\sigma,\, X,\, Y)\,:\;{\rm card\,}\sigma\leq n\,,\,\left|\lambda\right|\leq r,\;\forall\lambda\in\sigma\right\} .\]

\begin{flushleft}
{\small \vspace{0.1cm}
}
\par\end{flushleft}{\small \par}

It is explained in {[}15{]} why the classical interpolation problems,
those of Nevanlinna-Pick and Carathéodory-Schur (see {[}12{]} p.231),
on the one hand and Carleson's free interpolation (1958) (see {[}13{]}
p.158) on the other hand, are of this nature. 

From now on, if $\sigma=\left\{ \lambda_{1},\,...,\,\lambda_{n}\right\} \subset\mathbb{D}$
is a finite subset of the unit disc, then \[
B_{\sigma}={\displaystyle \prod_{j=1}^{n}}b_{\lambda_{j}}\]
 is the corresponding finite Blaschke product where $b_{\lambda}=\frac{\lambda-z}{1-\overline{\lambda}z}\:,$
$\lambda\in\mathbb{D}$. With this notation and supposing that $X$
satisfies the division property\[
\left[f\in X,\,\lambda\in\mathbb{D}\;{\rm and}\; f(\lambda)=0\right]\Rightarrow\left[\frac{f}{z-\lambda}\in X\right],\]
 we have 

\begin{onehalfspace}
\[
c\left(\sigma,\, X,\, Y\right)=\sup_{\parallel f\parallel_{X}\leq1}{\rm inf}\left\{ \left\Vert g\right\Vert _{Y}:\, g\in Y,\, g-f\in B_{\sigma}X\right\} .\]

\end{onehalfspace}

\subsection*{b. Motivations in matrix analysis and in operator theory}

A direct relation between the study of the constants $c\left(\sigma,\, H^{\infty},\, W\right)$
and some numerical analysis problems is mentioned in {[}15{]} (page
5, (b)). Here, $W$ is the Wiener algebra of absolutely convergent
Fourier series. In the same spirit, for general Banach spaces $X$
containing $H^{\infty}$, our constants $c\left(\sigma,\, X,\, H^{\infty}\right)$
are directly linked with the well known Von-Neumann's inequality for
contractions on Hilbert spaces, which asserts that if $A$ is a contraction
on a Hilbert space and $f\in H^{\infty},$ then the operator $f(A)$
satisfies\[
\left\Vert f(A)\right\Vert \leq\left\Vert f\right\Vert _{\infty}.\]
Using this inequality we get the following interpretation of our interpolation
constant $c\left(\sigma,\, X,\, H^{\infty}\right)$: it is the best
possible constant $c$ such that $\left\Vert f(A)\right\Vert \leq c\left\Vert f\right\Vert _{X},\;\forall\, f\in X$.
That is to say: 

\vspace{0.2cm}
\[
c\left(\sigma,\, X,\, H^{\infty}\right)=\sup_{\parallel f\parallel_{X}\leq1}\mbox{sup}\left\{ \left\Vert f(A)\right\Vert :\, A:\left(\mathbb{C}^{n},\,\vert\cdot\vert_{2}\right)\rightarrow\left(\mathbb{C}^{n},\,\vert\cdot\vert_{2}\right),\,\left\Vert A\right\Vert \leq1,\,\sigma(A)\subset\sigma\right\} ,\]

\begin{flushleft}
\vspace{0.5cm}
where the interior sup is taken over all contractions $A$ on $n-$dimensional
Hilbert spaces $\left(\mathbb{C}^{n},\,\vert.\vert_{2}\right)$, with
a given spectrum $\sigma(A)\subset\sigma$. 
\par\end{flushleft}

An interesting case occurs for $f$ such that $f_{\vert\sigma}=(1/z)_{\vert\sigma}$
(estimates on condition numbers and the norm of inverses of $n\times n$
matrices) or $f_{\vert\sigma}={\displaystyle [1/(\lambda-z)]_{\vert\sigma}}$
(estimates on the norm of the resolvent of an $n\times n$ matrix),
see for instance {[}18{]}.

\subsection*{c. Known results}

Let $H^{p}$ ($1\leq p\leq\infty$) be the standard Hardy spaces and
let $L_{a}^{2}$ be the Bergman space on $\mathbb{D}$. We obtained
in {[}16{]} some estimates on $c\left(\sigma,\, X,\, H^{\infty}\right)$
for the cases $X\in\left\{ H^{p},\, L_{a}^{2}\right\} $.

\begin{flushleft}
\textbf{Theorem A}\textit{. Let $n\geq1$, $r\in[0,\,1)$, $p\in[1,\,+\infty]$
and $\left|\lambda\right|\leq r$. Then}
\par\end{flushleft}

\def\theequation{${1}$}\begin{equation}
\frac{1}{32^{\frac{1}{p}}}\left(\frac{n}{1-\left|\lambda\right|}\right)^{\frac{1}{p}}\leq c\left(\sigma_{n,\,\lambda},\, H^{p},H^{\infty}\right)\leq C_{n,r}\left(H^{p},\, H^{\infty}\right)\leq A_{p}\left(\frac{n}{1-r}\right)^{\frac{1}{p}},\label{eq:}\end{equation}

\def\theequation{${2}$}\begin{equation}
\frac{1}{32}\frac{n}{1-\left|\lambda\right|}\leq c\left(\sigma_{n,\,\lambda},\, L_{a}^{2},\, H^{\infty}\right)\leq C_{n,\, r}\left(L_{a}^{2},\, H^{\infty}\right)\leq\sqrt{2}10^{\frac{1}{4}}\frac{n}{1-r},\label{eq:}\end{equation}
\textit{where \[
\sigma_{n,\,\lambda}=\{\lambda,\,...,\,\lambda\},\;(n\; times),\]
 is the one-point set of multiplicity $n$ corresponding to $\lambda,$
$A_{p}$ is a constant depending only on $p$ and the left-hand side
inequality in (1) is valid only for $p\in2\mathbb{Z}_{+}.$ For $p=2,$
we have $A_{2}=\sqrt{2}$. }

\begin{flushleft}
Note that this theorem was partially motivated by a question posed
in an applied situation in {[}5, 6{]}.
\par\end{flushleft}

Trying to generalize inequalities (1) and (2) for general Banach spaces
$X$ (of analytic functions of moderate growth in $\mathbb{D}$),
we formulate the following conjecture: $C_{n,\, r}\left(X,\, H^{\infty}\right)\leq a\varphi_{X}\left(1-\frac{1-r}{n}\right)$,
where $a$ is a constant depending on $X$ only and where $\varphi_{X}(t)$
stands for the norm of the evaluation functional $f\mapsto f(t)$
on the space $X$. The aim of this paper is to establish this conjecture
for some families of weighted Hardy and Bergman spaces.

\section*{2. Main results}

Here, we extend Theorem A to the case where $X$ is a weighted space
\[
l_{a}^{p}(\alpha)=\left\{ f={\displaystyle \sum_{k\geq0}\hat{f}(k)z^{k}:\,\Vert f\Vert^{p}=\sum_{k\geq0}\vert\hat{f}(k)\vert^{p}(k+1)^{p\alpha}<\infty}\right\} ,\;\alpha\leq0.\]
First, we study the special case $p=2$, $\alpha\leq0$. Then $l_{a}^{p}(\alpha)$
are the spaces of the functions $f=\sum_{k\geq0}\hat{f}(k)z^{k}$
satisfying \[
\sum_{k\geq0}\vert\hat{f}(k)\vert^{2}(k+1)^{2\alpha}<\infty.\]
Notice that $H^{2}=l_{a}^{2}(1).$ Let $\;\beta=-2\alpha-1>-1$. The
scale of weighted Bergman spaces of holomorphic functions \[
X=L_{a}^{2}\left(\beta\right)=L_{a}^{2}\left(\left(1-\left|z\right|^{2}\right)^{\beta}dA\right)=\left\{ f\in{\rm Hol}(\mathbb{D})\,:\;\int_{\mathbb{D}}\left|f(z)\right|^{2}\left(1-\left|z\right|^{2}\right)^{\beta}dA<\infty\right\} ,\]
gives the same spaces, with equivalence of the norms:\[
l_{a}^{2}\left(\alpha\right)=L_{a}^{2}\left(\beta\right).\]
In the case $\beta=0$ we have $L_{a}^{2}\left(0\right)=L_{a}^{2}$. 

We start with the following result.

\begin{flushleft}
\textbf{Theorem B.} \textit{Let $n\geq1$, $r\in[0,\,1),$ $\alpha\in(-\infty,\,0]$
and $\left|\lambda\right|\leq r$. Then}
\par\end{flushleft}

\textit{\[
B\left(\frac{n}{1-\left|\lambda\right|}\right)^{\frac{1-2\alpha}{2}}\leq c\left(\sigma_{n,\,\lambda},\, l_{a}^{2}\left(\alpha\right),\, H^{\infty}\right)\leq C_{n,\, r}\left(l_{a}^{2}\left(\alpha\right),\, H^{\infty}\right)\leq A\left(\frac{n}{1-r}\right)^{\frac{1-2\alpha}{2}}.\]
Equivalently, if $\beta\in(-1,\,+\infty)$ then \[
B'\left(\frac{n}{1-\left|\lambda\right|}\right)^{\frac{\beta+2}{2}}\leq c\left(\sigma_{n,\,\lambda},\, L_{a}^{2}\left(\beta\right),\, H^{\infty}\right)\leq C_{n,\, r}\left(L_{a}^{2}\left(\beta\right),\, H^{\infty}\right)\leq A'\left(\frac{n}{1-r}\right)^{\frac{\beta+2}{2}},\]
where $A$ and $B$ depend only on $\alpha$, $A'$ and $B'$ depend
only on $\beta$, and both of the two left-hand side inequalities
are valid only for $\alpha$ and $\beta$ satisfying $1-2\alpha\in\mathbb{N}$
and $\frac{\beta+1}{2}\in\mathbb{N}$.}

The right-hand side inequalities given in Theorem B are proved in
Section 4 whereas the left-hand side ones are proved in Section 5
.

\begin{flushleft}
\textit{Remark.} If $N=[1-2\alpha]$ is the integer part of $1-2\alpha$,
then Theorem B is valid with $B$ and $A$ such that $B\asymp\frac{1}{2^{3N}(2N)!}$
and $A\asymp N!(4N)^{N}$. In the same way, if $N'=\left[2+\beta\right]$
is the integer part of $2+\beta$, then Theorem B is valid with $B'$
and $A'$ such that $B'\asymp\frac{1}{2^{3N'}(2N')!}$ and $A'\asymp N'!(4N')^{N'}$.
(The notation $x\asymp y$ means that there exist numerical constants
$c_{1},\, c_{2}>0$ such that $c_{1}y\leq x\leq c_{2}y$).
\par\end{flushleft}

Next, we give an estimate for $C_{n,\, r}\left(X,\, H^{\infty}\right)$
in the scale of the spaces $X=l_{a}^{p}\left(\alpha\right)$, $\alpha\leq0$,
$1\leq p\leq+\infty$. We start with a result for $1\leq p\leq2$. 

\begin{flushleft}
\textbf{Theorem C. }\textit{Let $r\in[0,\,1),$ $n\geq1,$ $p\in[1,\,2]$,
and let $\alpha\leq0$. We have\[
Bn^{1-\alpha-\frac{1}{p}}\leq C_{n,\, r}\left(l_{a}^{p}\left(\alpha\right),\, H^{\infty}\right)\leq A\left(\frac{n}{1-r}\right)^{\frac{1-2\alpha}{2}},\]
}
\par\end{flushleft}

\begin{flushleft}
\textit{\vspace{0.5cm}
 where $A=A(\alpha,\, p)$ and $B=B(\alpha,\, p)$ are constants depending
only on $\alpha$ and $p$. }
\par\end{flushleft}

It is very likely that the bounds stated in Theorem C are not sharp.
The sharp one should be probably $\left(\frac{n}{1-r}\right)^{1-\alpha-\frac{1}{p}}$.
In the same way, for $2\leq p\leq\infty$, we give the following theorem,
in which we feel again that the upper bound $\left(\frac{n}{1-r}\right)^{\frac{3}{2}-\alpha-\frac{2}{p}}$
is not sharp. As before, the sharp one is probably $\left(\frac{n}{1-r}\right)^{1-\alpha-\frac{1}{p}}$.

\begin{flushleft}
\textbf{Theorem D.} \textit{Let $r\in[0,\,1),$ $n\geq1,$ $p\in[2,\,+\infty]$,
and let $\alpha\leq0$. We have}
\par\end{flushleft}

\textit{\[
B^{'}n^{1-\alpha-\frac{1}{p}}\leq C_{n,\, r}\left(l_{a}^{p}\left(\alpha\right),\, H^{\infty}\right)\leq A^{'}\left(\frac{n}{1-r}\right)^{\frac{3}{2}-\alpha-\frac{2}{p}},\]
}

\begin{flushleft}
\textit{\vspace{0.5cm}
 where $A^{'}$ and $B^{'}$ depend only on $\alpha$ and $p$. }
\par\end{flushleft}

\vspace{0.3cm}

Theorems B, C and D were already announced in the note {[}17{]}. Let
$\sigma$ be a finite set of $\mathbb{D},$ and let $f\in X.$ The
technical tools used in the proofs of the upper bounds for the interpolation
constants $c\left(\sigma,\, X,\, H^{\infty}\right)$ are: a linear
interpolation

\[
f\mapsto\sum_{k=1}^{n}\left\langle f,\, e_{k}\right\rangle e_{k},\]
where $\left\langle .,.\right\rangle $ means the Cauchy sesquilinear
form $\left\langle h,\, g\right\rangle =\sum_{k\geq0}\hat{h}(k)\overline{\hat{g}(k)},$
and $\left(e_{k}\right)_{1\leq k\leq n}$ is the explicitly known
Malmquist basis (see {[}13{]} p. 117) or Definition 1.1 below) of
the space $K_{B}=H^{2}\Theta BH^{2}$ where $B=B_{\sigma}$ (Subsection
3.1), a Bernstein-type inequality of Dyakonov (used by induction):
$\left\Vert f'\right\Vert _{p}\leq c_{p}\left\Vert B'\right\Vert _{\infty}\left\Vert f\right\Vert _{p},$
for a (rational) function $f$ in the star-invariant subspace $H^{p}\cap B\overline{zH^{p}}$
generated by a (finite) Blaschke product $B$, (Dyakonov {[}9, 10{]});
it is used in order to find an upper bound for $\left\Vert \sum_{k=1}^{n}\left\langle f,\, e_{k}\right\rangle e_{k}\right\Vert _{\infty}$
(in terms of $\left\Vert f\right\Vert _{X}$) (Subsection 3.2), and
finally (Subsection 3.3) the complex interpolation between Banach
spaces, (see {[}4{]} or {[}14{]} Theorem 1.9.3-(a), p.59).

The lower bound problem (for $C_{n,\, r}\left(X,\, H^{\infty}\right)$)
is treated by using the {}``worst'' interpolation $n-$tuple $\sigma=\sigma_{n,\,\lambda}=\{\lambda,\,...,\,\lambda\}$,
a one-point set of multiplicity $n$ (the Carathéodory-Schur type
interpolation). The {}``worst'' interpolation data comes from the
Dirichlet kernels $\sum_{k=0}^{n-1}z^{k}$ transplanted from the origin
to $\lambda.$ We note that the spaces $X=l_{a}^{p}(\alpha)$ satisfy
the condition $X\circ b_{\lambda}\subset X$ when $p=2$, whereas
this is not the case for $p\neq2$. That is why our problem of estimating
the interpolation constants is more difficult for $p\neq2$. 

\vspace{0.1cm}
The paper is organized as follows. In Section 3, we introduce the
three technical tools mentioned above. Section 4 is devoted to the
proof of the upper bounds of Theorems B, C and D. Finally, in Section
5, we prove the lower bounds of these theorems.

\section*{3. Preliminaries}

In this section, we develop the technical tools mentioned in Section
2, which are used later on to establish an upper bound for $c\left(\sigma,\, X,\, H^{\infty}\right)$.

\subsection*{3.1. Malmquist basis and orthogonal projection}

In Definitions 3.1.1, 3.1.2, 3.1.3 and in Remark 3.1.4 below, $\sigma=\left\{ \lambda_{1},\,...,\,\lambda_{n}\right\} $
is a sequence in the unit disc $\mathbb{D}$ and $B_{\sigma}$ is
the corresponding Blaschke product.

\begin{flushleft}
\textbf{Definition 3.1.1.}\textit{ Malmquist family. }For $k\in[1,\, n]$,
we set $f_{k}=\frac{1}{1-\overline{\lambda_{k}}z},$ and define the
family $\left(e_{k}\right)_{1\leq k\leq n}$, (which is known as Malmquist
basis, see {[}13, p.117{]}), by
\par\end{flushleft}

\def\theequation{${3.1.1}$}\begin{equation}
e_{1}=\frac{f_{1}}{\left\Vert f_{1}\right\Vert _{2}}\,\,\,\mbox{and}\,\,\, e_{k}=\left({\displaystyle \prod_{j=1}^{k-1}}b_{\lambda_{j}}\right)\frac{f_{k}}{\left\Vert f_{k}\right\Vert _{2}}\,,\label{eq:}\end{equation}
for $k\in[2,\, n]$; we have $\left\Vert f_{k}\right\Vert _{2}=\left(1-\vert\lambda_{k}\vert^{2}\right)^{-1/2}.$

\begin{flushleft}
\textbf{Definition 3.1.2.}\textit{ The model space $K_{B_{\sigma}}$.
}We define $K_{B_{\sigma}}$ to be the $n$-dimensional space:
\par\end{flushleft}

\def\theequation{${3.1.2}$}\begin{equation}
K_{B_{\sigma}}=\left(B_{\sigma}H^{2}\right)^{\perp}=H^{2}\Theta B_{\sigma}H^{2}.\label{eq:}\end{equation}

\begin{flushleft}
\textbf{Definition 3.1.3.}\textit{ The orthogonal projection ~$P_{B_{\sigma}}$on
$K_{B_{\sigma}}.$ }We define $P_{B_{\sigma}}$ to be the orthogonal
projection of $H^{2}$ on its $n$-dimensional subspace $K_{B_{\sigma}}.$
\par\end{flushleft}

\begin{flushleft}
\textbf{Remark 3.1.4.} The Malmquist family $\left(e_{k}\right)_{1\leq k\leq n}$
corresponding to $\sigma$ is an orthonormal basis of $K_{B_{\sigma}}.$
In particular,
\par\end{flushleft}

\def\theequation{${3.1.4}$}\begin{equation}
P_{B_{\sigma}}=\sum_{k=1}^{n}\left(\cdot,\, e_{k}\right)_{H^{2}}e_{k}\,,\label{eq:}\end{equation}

\begin{flushleft}
where $\left(.,\,.\right)_{H^{2}}$ means the scalar product on $H^{2}$.
\par\end{flushleft}

We now recall the following lemma already (partially) established
in {[}15, Lemma 3.1.5{]} which is useful in the proof of the upper
bound in Theorem C. 

\begin{flushleft}
\textbf{Lemma 3.1.5.}\textit{ Let $\sigma=\left\{ \lambda_{1},\,...,\,\lambda_{n}\right\} $
be a sequence in the unit disc $\mathbb{D}$ and let $\left(e_{k}\right)_{1\leq k\leq n}$
be the Malmquist family corresponding to $\sigma.$ Let also $\left\langle \cdot,\,\cdot\right\rangle $
be the Cauchy sesquilinear form $\left\langle h,\, g\right\rangle =\sum_{k\geq0}\hat{h}(k)\overline{\hat{g}(k)},$
(if $h\in{\rm Hol}(\mathbb{D})$ and $k\in\mathbb{N}$, $\hat{h}(k)$
stands for the $k^{th}$ Taylor coefficient of $h$). The map $\Phi:\,{\rm Hol}(\mathbb{D})\rightarrow{\rm Hol}(\mathbb{D})$
given by}
\par\end{flushleft}

\textit{\[
\:\:\:\:\:\:\:\Phi:\,\, f\mapsto\sum_{k=1}^{n}\left\langle f,\, e_{k}\right\rangle e_{k},\]
is well defined and has the following properties:}

\textit{(a) }\textbf{\textit{$\Phi_{\vert H^{2}}=P_{B_{\sigma}},$}}\textit{ }

\textit{(b) $\Phi$ is continuous on ${\rm Hol}(\mathbb{D})$ with
the topology of the uniform convergence on compact sets of $\mathbb{D}$,}

\begin{onehalfspace}
\textit{(c) if $X=l_{a}^{p}(\alpha)$ with $p\in[1,\,+\infty],$ $\alpha\in(-\infty,\,0]$
and $\Psi=Id_{\vert X}-\Phi_{\vert X},$ then ${\rm Im}\left(\Psi\right)\subset B_{\sigma}X,$}
\end{onehalfspace}

\textit{(d) if $f\in Hol(\mathbb{D}),$ then \[
\left|\Phi(f)(\zeta)\right|=\left|\left\langle f,\, P_{B_{\sigma}}k_{\zeta}\right\rangle \right|,\]
for all $\zeta\in\mathbb{D},$ where $P_{B_{\sigma}}$ is defined
in 3.1.3 and $k_{\zeta}=\left(1-\overline{\zeta}z\right)^{-1}.$}

\begin{flushleft}
\textit{Proof. }Points (a), (b) and (c) were already proved in {[}15{]}.
In order to prove (d), we simply need to write that \[
\Phi(f)(\zeta)=\sum_{k=1}^{n}\left\langle f,\, e_{k}\right\rangle e_{k}(\zeta)=\left\langle f,\,\sum_{k=1}^{n}\overline{e_{k}(\zeta)}e_{k}\right\rangle ,\]
\textit{$\forall f\in{\rm Hol}(\mathbb{D}),$ $\forall\zeta\in\mathbb{D}$
}and to notice that \textit{$\sum_{k=1}^{n}\overline{e_{k}(\zeta)}e_{k}=\sum_{k=1}^{n}\left(k_{\zeta},\, e_{k}\right)_{H^{2}}e_{k}=P_{B_{\sigma}}k_{\zeta}.$}
\par\end{flushleft}

\begin{flushright}
$\square$
\par\end{flushright}

\subsection*{3.2. Bernstein-type inequalities for rational functions}

Bernstein-type inequalities for rational functions are the subject
of a number of papers and monographs (see, for instance, {[}2, 3,
7, 8, 11{]}). We use here a result going back to Dyakonov {[}9, 10{]}.

\begin{flushleft}
\textbf{Lemma 3.2.1.} \textit{Let $B=\prod_{j=1}^{n}b_{\lambda_{j}}$,
be a finite Blaschke product (of order n), }$r=\max_{j}\left|\lambda_{j}\right|,$\textit{
and let $f\in K_{B}.$ Then}\[
\left\Vert f'\right\Vert _{H^{2}}\leq3\frac{n}{1-r}\left\Vert f\right\Vert _{H^{2}}.\]
Lemma 3.2.1 is a partial case ($p=2$) of the following K. Dyakonov's
result {[}8{]} (which is, in turn, a generalization of Levin's inequality
{[}11{]} corresponding to the case $p=\infty$): the norm $\left\Vert D\right\Vert _{K_{B}^{p}\rightarrow H^{p}}$
of the differentiation operator $Df=f^{'}$ on the star-invariant
subspace of the Hardy space $H^{p}$, $K_{B}^{p}:=H^{p}\cap B\overline{zH^{p}}$,
(where the bar denotes complex conjugation) satisfies the following
estimate: \[
\left\Vert D\right\Vert _{K_{B}^{p}\rightarrow H^{p}}\leq c_{p}\left\Vert B'\right\Vert _{\infty},\]
for every $p$, $1\leq p\leq\infty,$ where $c_{p}$ is a positive
constant depending only on $p$, $B$ is a finite Blaschke product
and $\left\Vert \cdot\right\Vert _{\infty}$ means the norm in $L^{\infty}(\mathbb{T})$.
In the case $p=2,$ Dyakonov's result gives $c_{p}=\frac{36+2\sqrt{3\pi}}{2\pi},$
which entails an estimate similar to that of Lemma 3.2.1, but with
a larger constant $\left(\frac{13}{2}\;{\rm instead\; of\;}3\right).$
Our lemma is proved in {[}16{]}, Proposition 6.1.1.
\par\end{flushleft}

The sharpness of the inequality stated in Lemma 3.2.1 is discussed
in {[}15{]}. Here we use it by induction in order to get the following
corollary.

\begin{flushleft}
\textbf{Corollary 3.2.2.} \textit{Let $B=\prod_{j=1}^{n}b_{\lambda_{j}}$,
be a finite Blaschke product (of order n), }$r=\max_{j}\left|\lambda_{j}\right|,$
\textit{and $f\in K_{B}.$ Then,}\[
\left\Vert f^{(k)}\right\Vert _{H^{2}}\leq k!4^{k}\left(\frac{n}{1-r}\right)^{k}\left\Vert f\right\Vert _{H^{2}},\]
\textit{for every $k=0,\,1,\,...$}
\par\end{flushleft}

\begin{flushleft}
\textit{Proof.} Indeed, since $z^{k-1}f^{(k-1)}\in K_{B^{k}},$ we
obtain applying Lemma 3.2.1 with $B^{k}$ instead of $B$,\[
\left\Vert z^{k-1}f^{(k)}+(k-1)z^{k-2}f^{(k-1)}\right\Vert _{H^{2}}\leq3\frac{kn}{1-r}\left\Vert z^{k-1}f^{(k-1)}\right\Vert _{H^{2}}=3\frac{kn}{1-r}\left\Vert f^{(k-1)}\right\Vert _{H^{2}}.\]
 In particular, \[
\left|\left\Vert z^{k-1}f^{(k)}\right\Vert _{H^{2}}-\left\Vert (k-1)z^{k-2}f^{(k-1)}\right\Vert _{H^{2}}\right|\leq3\frac{kn}{1-r}\left\Vert f^{(k-1)}\right\Vert _{H^{2}},\]
which gives \[
\left\Vert f^{(k)}\right\Vert _{H^{2}}\leq3\frac{kn}{1-r}\left\Vert f^{(k-1)}\right\Vert _{H^{2}}+(k-1)\left\Vert f^{(k-1)}\right\Vert _{H^{2}}\leq4\frac{kn}{1-r}\left\Vert f^{(k-1)}\right\Vert _{H^{2}}\,.\]
By induction,\[
\left\Vert f^{(k)}\right\Vert _{H^{2}}\leq k!\left(\frac{4n}{1-r}\right)^{k}\left\Vert f\right\Vert _{H^{2}}.\]

\par\end{flushleft}

\begin{flushright}
$\square$
\par\end{flushright}

\subsection*{3.3. Interpolation between Banach spaces (the complex method)}

In Section 4 we use the following lemma.

\begin{flushleft}
\textbf{Lemma 3.3.}\textit{ Let $X_{1}$ and $X_{2}$ be two Banach
spaces of holomorphic functions in the unit disc $\mathbb{D}$. Let
also $\theta\in[0,\,1]$ and $\left(X_{1},\, X_{2}\right)_{[\theta]}$
be the corresponding intermediate Banach space resulting from the
classical complex interpolation method applied between $X_{1}$ and
$X_{2}$, (we use the notation of {[}4, Chapter 4{]}). Then,\[
C_{n,\, r}\left(\left(X_{1},\, X_{2}\right)_{[\theta]},\, H^{\infty}\right)\leq C_{n,\, r}\left(X_{1},\, H^{\infty}\right)^{1-\theta}C_{n,\, r}\left(X_{2},\, H^{\infty}\right)^{\theta},\]
}
\par\end{flushleft}

\begin{flushleft}
\textit{for all $n\geq1,\, r\in[0,\,1).$}
\par\end{flushleft}

\begin{flushleft}
\textit{Proof.} Let $X$ be a Banach space of holomorphic functions
in the unit disc $\mathbb{D}$ and let $\sigma=\{\lambda_{1},\,\lambda_{2},\,...,\,\lambda_{n}\}\subset\mathbb{D}$
be a finite subset of the disc. Let $T\::\, X\longrightarrow H^{\infty}/B_{\sigma}H^{\infty}$
be the restriction map defined by \[
Tf=\left\{ g\in H^{\infty}:\: f-g\in B_{\sigma}X\right\} ,\]
for every $f\in X$. Then,
\par\end{flushleft}

\[
\left\Vert T\right\Vert _{X\rightarrow H^{\infty}/B_{\sigma}H^{\infty}}=c\left(\sigma,\, X,\, H^{\infty}\right).\]
Now, since $\left(X_{1},\, X_{2}\right)_{[\theta]}$ is an exact interpolation
space of exponent $\theta$ (see {[}4{]} or {[}14{]} Theorem 1.9.3-(a),
p.59), we can complete the proof. 

\begin{flushright}
$\square$
\par\end{flushright}

\section*{4. upper bounds for $C_{n,\, r}\left(X,\, H^{\infty}\right)$\textbf{\noun{\Large{} }}}

The aim of this section is to prove the upper bounds stated in Theorems
B, C, and D.

\subsection*{4.1. The case $X=l_{a}^{2}\left(\alpha\right),\,\alpha\leq0$}

We start with the following result.

\begin{flushleft}
\textbf{Corollary 4.1.1.}\textbf{\textit{ }}\textit{Let $N\geq0$
be an integer. Then, }
\par\end{flushleft}

\textit{\[
C_{n,\, r}\left(l_{a}^{2}\left(-N\right),\, H^{\infty}\right)\leq A\left(\frac{n}{1-r}\right)^{\frac{2N+1}{2}},\]
for all $r\in[0,\,1[,$ $n\geq1,$ where $A$ depends only on $N$
(of order }$N!(4N)^{N},$\textit{ see the proof below). }

\begin{flushleft}
\textit{Proof.} Indeed, let $X=l_{a}^{2}\left(-N\right)$, $\sigma$
a finite subset of $\mathbb{D}$ and $B=B_{\sigma}.$ If $f\in X,$
then using part (c) of Lemma 3.1.5, we get that $\Phi(f)_{\vert\sigma}=f_{\vert\sigma}.$
Now, denoting $X^{\star}$ the dual of $X$ with respect to the Cauchy
pairing $\left\langle \cdot,\,\cdot\right\rangle $ (defined in Lemma
3.1.5). Applying point (d) of the same lemma, we obtain $X^{\star}=l_{a}^{2}\left(N\right)$
and \[
\left|\Phi(f)(\zeta)\right|\leq\left\Vert f\right\Vert _{X}\left\Vert P_{B}k_{\zeta}\right\Vert _{X^{\star}}\leq\left\Vert f\right\Vert _{X}K_{N}\left(\left\Vert P_{B}k_{\zeta}\right\Vert _{H^{2}}^{2}+\left\Vert \left(P_{B}k_{\zeta}\right)^{(N)}\right\Vert _{H^{2}}^{2}\right)^{\frac{1}{2}},\]
for all $\zeta\in\mathbb{D},$ where \[
K_{N}=\mbox{max}\left\{ N^{N},\,\sup_{k\geq N}\frac{(k+1)^{N}}{k(k-1)...(k-N+1)}\right\} =\]
\[
=\mbox{max}\left\{ N^{N},\,\frac{(N+1)^{N}}{N!}\right\} =\left\{ \begin{array}{c}
N^{N},\: if\: N\geq3\\
\frac{(N+1)^{N}}{N!},\: if\: N=1,\,2\end{array}\right..\]
(Indeed, the sequence $\left(\frac{(k+1)^{N}}{k(k-1)...(k-N+1)}\right)_{k\geq N}$
is decreasing and $\left[N^{N}>\frac{(N+1)^{N}}{N!}\right]\Longleftrightarrow N\geq3$).
Since $P_{B}k_{\zeta}\in K_{B}$, Corollary 3.2.2 implies\[
\left|\Phi(f)(\zeta)\right|\leq\left\Vert f\right\Vert _{X}K_{N}\left\Vert P_{B}k_{\zeta}\right\Vert _{H^{2}}\left(1+(N!)^{2}\left(4\frac{n}{1-r}\right)^{2N}\right)^{\frac{1}{2}}\leq A(N)\left(\frac{n}{1-r}\right)^{N+\frac{1}{2}}\left\Vert f\right\Vert _{X},\]
where $A(N)=\sqrt{2}K_{N}\left(1+(N!)^{2}4^{2N}\right)^{\frac{1}{2}},$
since
\par\end{flushleft}

\def\theequation{${4.1.2}$}\begin{equation}
\left\Vert P_{B}k_{\zeta}\right\Vert _{H^{2}}=\left\Vert \sum_{k=1}^{n}\left(k_{\zeta},\, e_{k}\right)_{H^{2}}e_{k}\right\Vert _{H^{2}}=\sqrt{\sum_{k=1}^{n}\left|e_{k}(\zeta)\right|^{2}}\leq\sqrt{\frac{2n}{1-r}.}\label{eq:}\end{equation}

\begin{flushright}
$\square$
\par\end{flushright}

\vspace{0.2cm}

\begin{flushleft}
\textit{Proof of Theorem B (the right-hand side inequality). }There
exists{\large{} }an integer $N$ such that $N-1\leq-\alpha\leq N.$
In particular, there exists $0\leq\theta\leq1$ such that $\alpha=(1-\theta)(1-N)+\theta.(-N)$.
Since \[
\left(l_{a}^{2}\left(1-N\right),\, l_{a}^{2}\left(-N\right)\right)_{[\theta]}=l_{a}^{2}\left(\alpha\right),\]
(see {[}4, 14{]}), this gives, using Lemma 3.3 with $X_{1}=l_{a}^{2}\left(1-N\right)$
and $X_{2}=l_{a}^{2}\left(-N\right),$ and Corollary 4.1.1, that
\par\end{flushleft}

\[
C_{n,\, r}\left(l_{a}^{2}\left(\alpha\right),\, H^{\infty}\right)\leq A(N-1)^{1-\theta}A(N)^{\theta}\left(\frac{n}{1-r}\right)^{\frac{(2N-1)(1-\theta)}{2}+\frac{(2N+1)\theta}{2}}.\]
It remains to use that $\theta=1-\alpha-N$ and set $A(\alpha)=A(N-1)^{1-\theta}A(N)^{\theta}.$ 

\begin{flushright}
$\square$
\par\end{flushright}

\subsection*{4.2. An upper bound for $c\left(\sigma,\, l_{a}^{p}\left(\alpha\right),\, H^{\infty}\right),\:1\leq p\leq2$ }

The purpose of this subsection is to prove the right-hand side inequality
of Theorem C. We start with a partial case. 

\begin{flushleft}
\textbf{Lemma 4.2.1. }\textit{Let $N\geq0$ be an integer. Then }
\par\end{flushleft}

\textit{\[
C_{n,\, r}\left(l_{a}^{1}\left((-N\right),\, H^{\infty}\right)\leq A_{1}\left(\frac{n}{1-r}\right)^{N+\frac{1}{2}},\]
for all $r\in[0,\,1),$ $n\geq1,$ where $A_{1}$ depends only on
$N$ (it is of order }$N!(4N)^{N},$\textit{ see the proof below). }

\begin{flushleft}
\textit{Proof.} In fact, the proof is exactly the same as in Corollary
4.1.1: if $\sigma$ is a sequence in $\mathbb{D}$ with ${\rm card\,}\sigma\leq n,$
and $f\in l_{a}^{1}\left(-N\right)=X,$ then $X^{\star}=l_{a}^{\infty}\left(N\right)$
(the dual of $X$ with respect to the Cauchy pairing). Using Lemma
3.1.5 we still have $\Phi(f)_{\vert\sigma}=f_{\vert\sigma},$ and
for every $\zeta\in\mathbb{D},$\[
\left|\Phi(f)(\zeta)\right|\leq\left\Vert f\right\Vert _{X}\left\Vert P_{B}k_{\zeta}\right\Vert _{X^{\star}}\leq\]
\[
\leq\left\Vert f\right\Vert _{X}K_{N}\mbox{max}\left\{ \sup_{0\leq k\leq N-1}\left|\widehat{P_{B}k_{\zeta}}(k)\right|,\,\sup_{k\geq N}\left|\widehat{\left(P_{B}k_{\zeta}\right)^{(N)}}\left(k-N\right)\right|\right\} \leq\]
\[
\leq\left\Vert f\right\Vert _{X}K_{N}\,\mbox{max}\left\{ \left\Vert P_{B}k_{\zeta}\right\Vert _{H^{2}},\,\left\Vert \left(P_{B}k_{\zeta}\right)^{(N)}\right\Vert _{H^{2}}\right\} ,\]
where $K_{N}$ is defined in the the proof of Corollary 4.1.1. Since
$P_{B}k_{\zeta}\in K_{B}$, Corollary 3.2.2 implies that\[
\left|\Phi(f)(\zeta)\right|\leq\left\Vert f\right\Vert _{X}K_{N}\left\Vert P_{B}k_{\zeta}\right\Vert _{H^{2}}\left(1+N!4^{N}\left(\frac{n}{1-r}\right)^{N}\right),\]
for all $\zeta\in\mathbb{D},$ which completes the proof using (4.1.2)
and setting $A_{1}(N)=2\sqrt{2}N!4^{N}K_{N}$.
\par\end{flushleft}

\begin{flushright}
$\square$
\par\end{flushright}

\begin{flushleft}
\textit{Proof of Theorem C (the right-hand inequality). }
\par\end{flushleft}

\begin{flushleft}
\textbf{Step 1.} We start by proving the result for $p=1$ and for
all $\alpha\leq0.$ We use the same reasoning as in Theorem B except
that we replace $l_{a}^{2}(\alpha)$ by $l_{a}^{1}(\alpha).$ 
\par\end{flushleft}

\begin{flushleft}
\textbf{Step 2. }We now prove the result for $p\in[1,\,2]$ and for
all $\alpha\leq0:$ the scheme of this step is completely the same
as in Step 1, but we use this time the complex interpolation between
$l_{a}^{1}(\alpha)$ and $l_{a}^{2}(\alpha)$ (the classical Riesz-Thorin
Theorem {[}4, 14{]}). Applying Lemma 3.3 with $X_{1}=l_{a}^{1}\left(\alpha\right)$
and $X_{2}=l_{a}^{2}\left(\alpha\right)$, it suffices to\textbf{
}use Theorem B and Theorem C for the special case $p=1$ (already
proved in Step 1), to complete the proof of the right-hand side inequality.
\par\end{flushleft}

\begin{flushright}
$\square$
\par\end{flushright}

\subsection*{4.3. An upper bound for $c\left(\sigma,\, l_{a}^{p}\left(\alpha\right),\, H^{\infty}\right),\:2\leq p\leq+\infty$ }

Here, we prove the upper bound stated in Theorem D. As before, the
upper bound $\left(\frac{n}{1-r}\right)^{\frac{3}{2}-\alpha-\frac{2}{p}}$
is not as sharp as in Subsection 4.1. As in Subsection 4.2, we can
suppose the constant $\left(\frac{n}{1-r}\right)^{1-\alpha-\frac{1}{p}}$should
be again a sharp upper (and lower) bound for the quantity $C_{n,\, r}\left(l_{a}^{p}\left(\alpha\right),\, H^{\infty}\right),\:2\leq p\leq+\infty\,.$ 

First we prove the following partial case of Theorem D.

\begin{flushleft}
\textbf{Corollary 4.3.1. }\textit{Let $N\geq0$ be an integer. Then, }
\par\end{flushleft}

\textit{\[
C_{n,\, r}\left(l_{a}^{\infty}\left(-N\right),\, H^{\infty}\right)\leq A_{\infty}\left(\frac{n}{1-r}\right)^{N+\frac{3}{2}},\]
for all $r\in[0,\,1[,$ $n\geq1$, where $A_{\infty}$ depends only
on $N$ (it is of order $N!(4N)^{N},$ see the proof below).} 
\begin{proof}
We use literally the same method as in Corollary 4.1.1 and Lemma 4.2.1.
Indeed, if $\sigma=\left\{ \lambda_{1},\,...,\,\lambda_{n}\right\} $
is a sequence in the unit disc $\mathbb{D}$ and $f\in l_{a}^{\infty}\left(-N\right)=X,$
then $X^{\star}=l_{a}^{1}\left(N\right)$ and applying again Lemma
3.1.5 we get $\Phi(f)_{\vert\sigma}=f_{\vert\sigma}.$ For every $\zeta\in\mathbb{D},$
we have \[
\left|\Phi(f)(\zeta)\right|\leq\left\Vert f\right\Vert _{X}\left\Vert P_{B}k_{\zeta}\right\Vert _{X^{\star}}\leq\left\Vert f\right\Vert _{X}K_{N}\left(\left\Vert P_{B}k_{\zeta}\right\Vert _{W}+\left\Vert \left(P_{B}k_{\zeta}\right)^{(N)}\right\Vert _{W}\right),\]
where $W=\left\{ f=\sum_{k\geq0}\hat{f}(k)z^{k}:\,\left\Vert f\right\Vert _{W}:=\sum_{k\geq0}\left|\hat{f}(k)\right|<\infty\right\} $
stands for the Wiener algebra, and $K_{N}$ is defined in the proof
of Corollary 4.1.1. Now, applying Hardy's inequality (see {[}13{]},
p.370{]}), we obtain\[
\left|\Phi(f)(\zeta)\right|\leq\]
\[
\leq\left\Vert f\right\Vert _{X}K_{N}\left(\pi\left\Vert \left(P_{B}k_{\zeta}\right)^{'}\right\Vert _{H^{1}}+\left|\left(P_{B}k_{\zeta}\right)(0)\right|+\pi\left\Vert \left(P_{B}k_{\zeta}\right)^{(N+1)}\right\Vert _{H^{1}}+\left|\left(P_{B}k_{\zeta}\right)^{(N)}(0)\right|\right)\leq\]
\[
\leq\left\Vert f\right\Vert _{X}K_{N}\pi\left(\left\Vert \left(P_{B}k_{\zeta}\right)^{'}\right\Vert _{H^{2}}+\left\Vert \left(P_{B}k_{\zeta}\right)\right\Vert _{H^{2}}+\left\Vert \left(P_{B}k_{\zeta}\right)^{(N+1)}\right\Vert _{H^{2}}+\left\Vert \left(P_{B}k_{\zeta}\right)^{(N)}\right\Vert _{H^{2}}\right),\]
for all $\zeta\in\mathbb{D}.$ Using Lemma 3.2.1 and Corollary 3.2.2,
we get\[
\left|\Phi(f)(\zeta)\right|\leq\left\Vert f\right\Vert _{X}K_{N}\pi\left\Vert P_{B}k_{\zeta}\right\Vert _{H^{2}}\left(\frac{3n}{1-r}+1+(N+1)!\left(\frac{4n}{1-r}\right)^{N+1}+N!\right),\]
for all $\zeta\in\mathbb{D},$ which completes the proof using (4.1.2).
\end{proof}
\begin{flushleft}
\textit{Proof of Theorem D (the right-hand side inequality).} The
proof repeates the scheme from Theorem C (the two steps) excepted
that this time, we replace (in both steps) the space $X=l_{a}^{1}(\alpha)$
by $X=l_{a}^{\infty}(\alpha).$ 
\par\end{flushleft}

\begin{flushright}
$\square$
\par\end{flushright}

\section*{5. Lower bounds for $C_{n,\, r}\left(X,\, H^{\infty}\right)$}

Here we prove the left-hand side inequalities stated in Theorems B,
C and D.

\subsection*{5.1. The case $X=l_{a}^{2}(\alpha),\,\alpha\leq0$}

We start with verifying the sharpness of the upper estimate for the
quantity

\[
C_{n,\, r}\left(l_{a}^{2}\left(\frac{1-N}{2}\right),\, H^{\infty}\right),\]
(where $N\geq1$ is an integer), in Theorem B. This lower bound problem
is treated by estimating our interpolation constant $c(\sigma,\; X,\; H^{\infty})$
for the one-point interpolation set $\sigma_{n,\,\lambda}=\underbrace{\{\lambda,\lambda,...,\lambda\}}_{n}$,
$\lambda\in\mathbb{D}$:

\[
c(\sigma_{n,\,\lambda},\; X,\; H^{\infty})=\mbox{sup}\left\{ \left\Vert f\right\Vert _{H^{\infty}/b_{\lambda}^{n}H^{\infty}}:\, f\in X,\,\left\Vert f\right\Vert _{X}\leq1\right\} ,\]

\begin{flushleft}
\vspace{0.3cm}
where $\left\Vert f\right\Vert _{H^{\infty}/b_{\lambda}^{n}H^{\infty}}=\mbox{inf}\left\{ \left\Vert f+b_{\lambda}^{n}g\right\Vert _{\infty}\,:\, g\in X\right\} $.
In the proof, we notice that $l_{a}^{2}(\alpha)$ is a reproducing
kernel Hilbert space on the disc $\mathbb{D}$ (RKHS) and we use the
fact that this space has some special properties for particular values
of $\alpha$ $\left(\alpha=\frac{1-N}{2},\; N=1,\,2,\,...\right)$.
Before giving this proof (see Paragraph 5.1.2 below), we show in Subsection
5.1.1 that $l_{a}^{2}(\alpha)$ is a RKHS and we focus on the special
case $\alpha=\frac{1-N}{2},\; N=1,\,2,\,...$.
\par\end{flushleft}

\subsubsection*{5.1.1. The spaces $l_{a}^{2}(\alpha)$ are RKHS}

The reproducing kernel of $l_{a}^{2}(\alpha)$, by definition, is
a $l_{a}^{2}(\alpha)$-valued function $\lambda\longmapsto k_{\lambda}^{\alpha}$,
$\lambda\in\mathbb{D}$, such that $\left(f,\, k_{\lambda}^{\alpha}\right)=f(\lambda)$
for every $f\in l_{a}^{2}(\alpha)$, where $\left(.,.\right)$ means
the scalar product $\left(f,\, g\right)=\sum_{k\geq0}\hat{h}(k)\overline{\hat{g}(k)}(k+1)^{2\alpha}.$
Since one has $f(\lambda)=\sum_{k\geq0}\hat{f}(k)\lambda^{k}\frac{1}{(k+1)^{2\alpha}}(k+1)^{2\alpha}$
($\lambda\in\mathbb{D}$), it follows that

\[
k_{\lambda}^{\alpha}(z)={\displaystyle \sum_{k\geq0}\frac{{\displaystyle \overline{\lambda}^{k}z^{k}}}{(k+1)^{2\alpha}}},\: z\in\mathbb{D}.\]
In particular, for the Hardy space $H^{2}=\, l_{a}^{2}(1)$, we get
the Szegö kernel\[
k_{\lambda}(z)=\,(1-\overline{\lambda}z)^{-1},\]
and for the Bergman space $L_{a}^{2}=\, l_{a}^{2}\left(-\frac{1}{2}\right),$
the Bergman kernel $k_{\lambda}^{-1/2}(z)=$ $(1-\overline{\lambda}z)^{-2}$.

Now let us explain that more generally if $\alpha=\frac{1-N}{2},\; N\in\mathbb{N}\setminus\{0\},$
the space $l_{a}^{2}(\alpha)$ coincides (topologically) with the
RKHS whose reproducing kernel is $\left(k_{\lambda}(z)\right)^{N}=(1-\overline{\lambda}z)^{-N}.$
Following the Aronszajn theory of RKHS (see, for example {[}1, 12{]}),
given a positive definite function $(\lambda,z)\longmapsto k(\lambda,\, z)$
on $\mathbb{D}\times\mathbb{D}$ (i.e. such that $\sum_{i,j}\overline{a}_{i}a_{j}k(\lambda_{i},\lambda_{j})>0$
for all finite subsets $(\lambda_{i})\subset\mathbb{D}$ and all non-zero
families of complex numbers $(a_{i})$) one can define the corresponding
Hilbert spaces $H(k)$ as the completion of finite linear combinations
$\sum_{i}\overline{a}_{i}k(\lambda_{i},\cdot)$ endowed with the norm\[
{\displaystyle \left\Vert \sum_{i}\overline{a}_{i}k(\lambda_{i},\cdot)\right\Vert ^{2}={\displaystyle \sum_{i,j}\overline{a}_{i}a_{j}k(\lambda_{i},\lambda_{j}).}}\]
When $k$ is holomorphic with respect to the second variable and antiholomorphic
with respect to the first one, we obtain a RKHS of holomorphic functions
$H(k)$ embedded into ${\rm Hol}(\mathbb{D})$. Now, choosing for
$k$ the reproducing kernel of $H^{2},$ $k\,:(\lambda,\, z)\mapsto k_{\lambda}(z)=\,(1-\overline{\lambda}z)^{-1},$
and $\varphi=z^{N},\; N=1,\,2,\,...,$ the function $\varphi\circ k$
is also positive definite and the corresponding Hilbert space is

\def\theequation{${5.1.1}$}\begin{equation}
H(\varphi\circ k)=l_{a}^{2}\left(\frac{1-N}{2}\right).\label{eq:}\end{equation}
(Another notation for the space $H(\varphi\circ k)$ is $\varphi(H^{2})$
since $k$ is the reproducing kernel of $H^{2}$). The equality (5.1.1)
is a topological identity: the spaces coincide as sets of functions,
and the norms are equivalent. Moreover, the space $H(\varphi\circ k)$
satisfies the following property: for every $f\in H^{2}$, $\varphi\circ f\in\varphi(H^{2}),$
and

\def\theequation{${5.1.2}$}\begin{equation}
\Vert\varphi\circ f\Vert_{H(\varphi\circ k)}^{2}\leq\varphi(\Vert f\Vert_{H^{2}}^{2}),\label{eq:}\end{equation}
(the Aronszajn-deBranges inequality, see {[}13{]} p.320{]}). The link
between spaces of type $l_{a}^{2}\left(\frac{1-N}{2}\right)$ and
of type $H(z^{N}\circ k)$ being established, we give the proof of
the left-hand side inequality in Theorem B.

\subsubsection*{5.1.2. The proof of Theorem B (the lower bound)}

0) We set $N=1-2\alpha,\; N=1,\,2,\,...$ and $\varphi(z)=z^{N}.$

1) Let $b>0,$ $b^{2}n^{N}=1.$ We set\[
Q_{n}=\sum_{k=0}^{n-1}b_{\lambda}^{k}\frac{(1-\vert\lambda\vert^{2})^{1/2}}{1-\overline{\lambda}z},\: H_{n}=\varphi\circ Q_{n},\;\Psi=bH_{n}.\]
Then $\Vert Q_{n}\Vert_{2}^{2}=\, n$, and hence by (5.1.2), 

\[
\Vert\Psi\Vert_{H_{\varphi}}^{2}\leq b^{2}\varphi\left(\Vert Q_{n}\Vert_{2}^{2}\right)=b^{2}\varphi(n)=1.\]

\vspace{1cm}
 Let $b>0$ such that $b^{2}\varphi(n)=\,1$.

\vspace{1cm}
 2) Since the spaces $H_{\varphi}$ and $H^{\infty}$ are rotation
invariant, we have $c\left(\sigma_{n,\,\lambda},\, H_{\varphi},\, H^{\infty}\right)=c\left(\sigma_{n,\,\mu},\, H_{\varphi},\, H^{\infty}\right)$
for every $\lambda,\,\mu$ with $\vert\lambda\vert=\vert\mu\vert=r$.
Let $\lambda=-r$. To get a lower estimate for $\Vert\Psi\Vert_{H_{\varphi}/b_{\lambda}^{n}H_{\varphi}}$
consider $G\in H^{\infty}$ such that $\Psi-G\in b_{\lambda}^{n}{\rm Hol}(\mathbb{D})$,
i.e. such that $bH_{n}\circ b_{\lambda}-G\circ b_{\lambda}\in z^{n}{\rm Hol}(\mathbb{D})$.

\vspace{1cm}
 3) First, we show that

\[
\psi=:\,\Psi\circ b_{\lambda}=bH_{n}\circ b_{\lambda}\]
is a polynomial (of degree $nN$) with positive coefficients. Note
that

\[
Q_{n}\circ b_{\lambda}=\sum_{k=0}^{n-1}z^{k}\frac{(1-\vert\lambda\vert^{2})^{1/2}}{1-\overline{\lambda}b_{\lambda}(z)}=\]

\[
=\left(1-\vert\lambda\vert^{2}\right)^{-\frac{1}{2}}\left(1+(1-\overline{\lambda})\sum_{k=1}^{n-1}z^{k}-\overline{\lambda}z^{n}\right)=\]
\[
=(1-r^{2})^{-1/2}\left(1+(1+r)\sum_{k=1}^{n-1}z^{k}+rz^{n}\right)=:(1-r^{2})^{-1/2}\psi_{1}.\]
Then, $\psi=\Psi\circ b_{\lambda}=bH_{n}\circ b_{\lambda}=b\varphi\circ\left(\left(1-r^{2}\right)^{-\frac{1}{2}}\psi_{1}\right).$
Furthermore, 

\[
\varphi\circ\psi_{1}=\psi_{1}^{N}(z).\]
Now, it is clear that $\psi$ is a polynomial of degree $Nn$ such
that\[
\psi(1)=\sum_{j=0}^{Nn}\hat{\psi}(j)=b\varphi\left((1-r^{2})^{-1/2}(1+r)n\right)=b\left(\sqrt{\frac{1+r}{1-r}}n\right)^{N}>0.\]

\vspace{1cm}
 4) Next, we show that there exists $c=c(N)>0$ (for example, $c=K/\left[2^{2N}(N-1)!\right]$,
$K$ being a numerical constant) such that\[
\sum^{m}(\psi)\,:=\sum_{j=0}^{m}\hat{\psi}(j)\geq c\sum_{j=0}^{Nn}\hat{\psi}(j)=c\psi(1),\]
where $m\geq1$ is such that $2m=n$ if $n$ is even and $2m-1=n$
if $n$ is odd. 

\vspace{1cm}
Indeed, setting\[
S_{n}=\sum_{j=0}^{n}z^{j},\]
we have\[
\sum^{m}\left(\psi_{1}^{N}\right)=\sum^{m}\left(\left(1+(1+r)\sum_{k=1}^{n-1}z^{k}+rz^{n}\right)^{N}\right)\geq\sum^{m}\left(S_{n-1}^{N}\right).\]
Next, we obtain\[
\sum^{m}\left(S_{n-1}^{N}\right)=\sum^{m}\left(\left(\frac{1-z^{n}}{1-z}\right)^{N}\right)=\]
\[
=\sum^{m}\left(\frac{1}{(1-z)^{N}}\right)=\frac{1}{(N-1)!}\sum^{m}\left(\frac{d^{N-1}}{dz^{N-1}}\frac{1}{1-z}\right)=\]
\[
=\sum_{j=0}^{m}C_{N+j-1}^{j}\geq\sum_{j=0}^{m}\frac{(j+1)^{N-1}}{(N-1)!}\geq\]
\[
\geq K\frac{m^{N}}{(N-1)!},\]
where $K>0$ is a numerical constant. Finally,\[
\sum^{m}\left(\psi_{1}^{N}\right)\geq K\frac{m^{N}}{(N-1)!}\geq K\frac{(n/2)^{N}}{(N-1)!}=\]
\[
=\frac{K}{2^{N}(N-1)!}\cdot\frac{((1+r)n)^{N}}{(1+r)^{N}}=\frac{K}{2^{N}(1+r)^{N}(N-1)!}\cdot(\psi_{1}(1))^{N},\]
which gives our estimate.

\vspace{1cm}

5) Let $F_{n}=\Phi_{m}+z^{m}\Phi_{m}$, where $\Phi_{k}$ stands for
the $k$-th Fejer kernel. We have $\Vert g\Vert_{\infty}\Vert F_{n}\Vert_{L^{1}}\geq\Vert g\star F_{n}\Vert_{\infty}$
for every $g\in L^{\infty}\left(\mathbb{T}\right),$ and taking the
infimum over all $g\in H^{\infty}$ satisfying $\hat{g}(k)=\hat{\psi}(k),\;\forall k\in[0,\, n-1],$
we obtain\[
\Vert\psi\Vert_{H^{\infty}/z^{n}H^{\infty}}\geq\frac{{\displaystyle 1}}{{\displaystyle 2}}\Vert\psi\star F_{n}\Vert_{\infty},\]
where $\star$ stands for the usual convolution product. Now using
part 4),

\[
\Vert\Psi\Vert_{H^{\infty}/b_{\lambda}^{n}H^{\infty}}=\Vert\psi\Vert_{H^{\infty}/z^{n}H^{\infty}}\geq\frac{{\displaystyle 1}}{{\displaystyle 2}}\Vert\psi\star F_{n}\Vert_{\infty}\geq\]

\[
\geq\frac{{\displaystyle 1}}{{\displaystyle 2}}\left|\left(\psi\star F_{n}\right)(1)\right|\geq\frac{{\displaystyle 1}}{{\displaystyle 2}}{\displaystyle \sum_{j=0}^{m}\hat{\psi}(j)}\geq\frac{{\displaystyle c}}{{\displaystyle 2}}\psi(1)=\frac{{\displaystyle c}}{{\displaystyle 2}}b\left(\sqrt{\frac{1+r}{1-r}}n\right)^{N}\geq\]

\[
\geq B\left(\frac{{\displaystyle n}}{{\displaystyle 1-r}}\right)^{\frac{N}{2}}.\]

6) In order to conclude, it remains to use (5.1.1).

\begin{flushright}
$\square$
\par\end{flushright}

\subsection*{5.2. The case $X=l_{a}^{p}(\alpha),\,1\leq p\leq\infty$ }

\subsection*{$\;$}

\begin{flushleft}
\textit{Proof of Theorems C and D (the lower bound)} We first notice
that $r\mapsto C_{n,\, r}\left(X,\, H^{\infty}\right)$ increases.
As a consequence, if $X=l_{a}^{p}(\alpha),\,1\leq p\leq\infty,$ then
\[
C_{n,\, r}\left(l_{a}^{p}\left(\alpha\right),\, H^{\infty}\right)\geq C_{n,\,0}\left(l_{a}^{p}\left(\alpha\right),\, H^{\infty}\right)=c\left(\sigma_{n,\,0},\, l_{a}^{p}\left(\alpha\right),\, H^{\infty}\right),\]
where $\sigma_{n,\,0}=\underbrace{\{0,\,0,...,\,0\}}_{n}.$ Now let
$f=\frac{1}{n^{1/p}}\sum_{k=0}^{n-1}(k+1)^{-\alpha}z^{k}$. Then $\left\Vert f\right\Vert _{X}=1,$
and \[
c\left(\sigma_{n,\,0},\, l_{a}^{p}\left(\alpha\right),\, H^{\infty}\right)\geq\Vert f\Vert_{H^{\infty}/z^{n}H^{\infty}}\geq\]
\[
\geq\frac{{\displaystyle 1}}{{\displaystyle 2}}\Vert f\star F_{n}\Vert_{\infty}\geq\frac{{\displaystyle 1}}{{\displaystyle 2}}\left|\left(f\star F_{n}\right)(1)\right|\geq\frac{{\displaystyle 1}}{{\displaystyle 2}}{\displaystyle \sum_{j=0}^{m}\hat{f}(j),}\]
where $\star$ and $F_{n}$ are defined in part 5) of the proof of
Theorem B (lower bound) in Subsection 5.1 and where $m\geq1$ is such
that $2m=n$ if $n$ is even and $2m-1=n$ if $n$ is odd as in part
4) of the proof of the same Theorem. Now, since\[
\sum_{j=0}^{m}\hat{f}(j)=\frac{1}{n^{1/p}}\sum_{k=0}^{m}(k+1)^{-\alpha},\]
we get the result.
\par\end{flushleft}

\begin{flushright}
$\square$
\par\end{flushright}

\thanks{\textbf{Acknowlgement.} I would like to thank Professor Nikolai Nikolski
for all of his work, his wisdom and the pleasure that our discussions
gave to me. I also would like to thank the referee for the careful
review and the valuable comments, which provided insights that helped
improve the paper. I am also deeply grateful to Professor Alexander
Borichev for the thorough, constructive and helpful remarks and suggestions
on the manuscript. }

\vspace{0.4cm}

\noun{CMI-LATP, UMR 6632, Université de Provence, 39, rue F.-Joliot-Curie,
13453 Marseille cedex 13, France}

\textit{E-mail address} : rzarouf@cmi.univ-mrs.fr

\begin{thebibliography}{18}
{\normalsize \bibitem[1]{key-1} N. Aronszajn, }\textit{\normalsize Theory
of reproducing kernels, }{\normalsize Transactions of American Mathematical
Society, 68:337-404, 1950.}{\normalsize \par}

{\normalsize \bibitem[2]{key-3} A. Baranov, }\textit{\normalsize Bernstein-type
inequalities for shift-coinvariant subspaces and their applications
to Carleson embeddings}{\normalsize . Journal of Functional Analysis
(2005) 223 (1): 116-146.}{\normalsize \par}

{\normalsize \bibitem[3]{key-4} A. Baranov, }\textit{\normalsize Embeddings
of model subspaces of the Hardy space: compactness and Schatten\textendash{}von
Neumann ideals}{\normalsize , Izv. Ross. Nauk Ser. Mat., translated
in Izv. Math. 73 (2009), no. 6, 1077-1100. }{\normalsize \par}

{\normalsize \bibitem[4]{key-3-3} J. Bergh , J. Löfström, }\textit{\normalsize Interpolation
Spaces. An Introduction}{\normalsize , Springer-Verlag (1976). }{\normalsize \par}

{\normalsize \bibitem[5]{key-1} L. Baratchart, }\textit{\normalsize Rational
and meromorphic approximation in Lp of the circle : system-theoretic
motivations, critical points and error rates.}{\normalsize{} In N.
Papamichael, S. Ruscheweyh, and E. Saff, editors, Computational Methods
and Function Theory, pages 45\textendash{}78. World Scientific Publish.
Co, 1999. }{\normalsize \par}

{\normalsize \bibitem[6]{key-2-1} L. Baratchart, F. Wielonsky, }\textit{\normalsize Rational
approximation problem in the real Hardy space $H_{2}$ and Stieltjes
integrals: a uniqueness theorem}{\normalsize , Constr. Approx. 9 (1993),
1-21.}{\normalsize \par}

{\normalsize \bibitem[7]{key-3-2} P. Borwein and T. Erdélyi, }\textit{\normalsize Polynomials
and Polynomial Inequalities}{\normalsize , Springer, New York, 1995. }{\normalsize \par}

{\normalsize \bibitem[8]{key-4-1} R. A. DeVore and G. G. Lorentz,}\textit{\normalsize{}
Constructive Approximation}{\normalsize , Springer-Verlag, Berlin,
1993.}{\normalsize \par}

{\normalsize \bibitem[9]{key-3-1} K. M. Dyakonov, }\textit{\normalsize Differentiation
in Star-Invariant Subspaces I. Boundedness and Compactness,}{\normalsize{}
J.Funct.Analysis, 192, 364-386, 2002. }{\normalsize \par}

{\normalsize \bibitem[10]{key-1-1-1} K. M. Dyakonov, }\textit{\normalsize Entire
functions of exponential type and model subspaces in $H^{p}$}{\normalsize ,
Zap. Nauchn. Sem. Leningrad. Otdel. Mat. Inst. Steklov. (LOMI) 190
(1991), 81-100 (Russian); translation in J. Math. Sci. 71, 2222-2233,
1994.}{\normalsize \par}

{\normalsize \bibitem[11]{key-1} M. B. Levin, }\textit{\normalsize Estimation
of the derivative of a meromorphic function on the boundary of the
domain}{\normalsize{} (Russian), Teor. Funkci\u{\i} Funkcional. Anal.
i Priložen. Vyp. 24, 68-85, 1975.}{\normalsize \par}

{\normalsize \bibitem[12]{key-10} N.Nikolski, }\textit{\normalsize Operators,
Function, and Systems : an easy reading}{\normalsize , Vol.1, Amer.
Math. Soc. Monographs and Surveys, 2002.}{\normalsize \par}

{\normalsize \bibitem[13]{key-9} N.Nikolski, }\textit{\normalsize Treatise
on the shift operator}{\normalsize , Springer-Verlag, Berlin etc.,
1986 (Transl. from Russian, }\textit{\normalsize Lekzii ob}{\normalsize{}
}\textit{\normalsize operatore sdviga}{\normalsize , {}``Nauja'',
Moskva, 1980).}{\normalsize \par}

{\normalsize \bibitem[14]{key-11} H. Triebel, }\textit{\normalsize Interpolation
theory, functions spaces, differential operators, }{\normalsize North-Holland
Publishing Comp., 1978.}{\normalsize \par}

{\normalsize \bibitem[15]{key-1-1} R. Zarouf, }\textit{\normalsize Asymptotic
sharpness of a Bernstein-type inequality for rational functions in
$H^{2}$,}{\normalsize{} to appear in St Petersburg Math. J.}{\normalsize \par}

{\normalsize \bibitem[16]{key-5} R. Zarouf, }\textit{\normalsize Effective
$H^{\infty}$ interpolation constrained by Hardy and Bergman norms,
}{\normalsize submitted. }{\normalsize \par}

{\normalsize \bibitem[17]{key-1} R. Zarouf, }\textit{\normalsize Interpolation
avec contraintes sur des ensembles finis du disque, }{\normalsize C.
R. Acad. Sci. Paris, Ser. I 347, 2009.}{\normalsize \par}

{\normalsize \bibitem[18]{key-1-4} R. Zarouf,}\textit{\normalsize{}
Sharpening a result by E.B. Davies and B. Simon}{\normalsize , C.
R. Acad. Sci. Paris, Ser. I 347 (2009).}{\normalsize \par}

\end{thebibliography}
\end{document}